\magnification=\magstep1
\hsize16truecm
\vsize23.5truecm
\topskip=1truecm
\raggedbottom
\abovedisplayskip=3mm
\belowdisplayskip=3mm
\abovedisplayshortskip=0mm
\belowdisplayshortskip=2mm
\normalbaselineskip=12pt
\normalbaselines
\font\titlefont= cmcsc10 at 12pt
\def\F{\Bbb F}
\def\R{\Bbb R}
\def\C{\Bbb C}
\def\Z{\Bbb Z}

\def\P{\Bbb P}

\def\vs{\vskip 1.0pc}
\def\ni{\noindent}

%
%
\catcode`\@=11
\font\tenmsa=msam10
\font\sevenmsa=msam7
\font\fivemsa=msam5
\font\tenmsb=msbm10
\font\sevenmsb=msbm7
\font\fivemsb=msbm5
\newfam\msafam
\newfam\msbfam
\textfont\msafam=\tenmsa  \scriptfont\msafam=\sevenmsa
  \scriptscriptfont\msafam=\fivemsa
\textfont\msbfam=\tenmsb  \scriptfont\msbfam=\sevenmsb
  \scriptscriptfont\msbfam=\fivemsb
\def\hexnumber@#1{\ifcase#1 0\or1\or2\or3\or4\or5\or6\or7\or8\or9\or
	A\or B\or C\or D\or E\or F\fi }
\edef\msa@{\hexnumber@\msafam}
\edef\msb@{\hexnumber@\msbfam}
\mathchardef\square="0\msa@03
\mathchardef\subsetneq="3\msb@28
\mathchardef\ltimes="2\msb@6E
\mathchardef\rtimes="2\msb@6F
\def\Bbb{\ifmmode\let\next\Bbb@\else
	\def\next{\errmessage{Use \string\Bbb\space only in math mode}}\fi\next}
\def\Bbb@#1{{\Bbb@@{#1}}}
\def\Bbb@@#1{\fam\msbfam#1}
\catcode`\@=12
%
%
%
%
\vskip 6.5pc
\noindent
\font\eighteenbf=cmbx10 scaled\magstep3
\vskip 2.0pc
\centerline{\eighteenbf  How to Construct   Curves Over Finite Fields}
\bigskip
\centerline{\eighteenbf With Many Points}
\noindent
\vskip 2pc
\font\titlefont=cmcsc10 at 12pt
\centerline{\titlefont Gerard van der Geer and Marcel van der Vlugt}
\vskip 2.0pc
\noindent
\bigskip
\noindent
\centerline{\bf  Introduction}
\smallskip

In 1940  A.\ Weil proved the Riemann hypothesis for curves over
finite fields. As an immediate corollary he obtained an upper bound for
the number of points on an irreducible  curve $C$ of genus g over a
finite field  of cardinality $q$, namely
$$
\# C(\F_q)  \leq q+1 + 2g\sqrt{q}.
$$
This bound was proved for elliptic curves by Hasse in 1933. Ever since,
the question of the maximum number
$N_q(g)$  of points on an irreducible  curve of genus $g$ over a
finite field of cardinality $q$ could have been investigated. But for
a long time it attracted no attention and it was only  after Goppa
introduced geometric codes in 1980 that this question  aroused
systematic attention, cf.\  [G], [M], [S 1,2,3].

For $g \leq (q- \sqrt{q})/2$ the Hasse-Weil bound is in general best
possible; for other $g$ there is a better bound. In general the
so-called `formules explicites' can provide improvements of the
Hasse-Weil bound. Therefore, the question arises of the  actual value
of the maximum number $N_q(g)$ for given $g$ and $q$.

Serre gave tables in [S 1,2,3] listing the value of $N_q(g)$ for small
values of $q$ and $g$, or if this value was not known, a small interval
in which the value of $N_q(g)$ lies. Wirtz extended in [Wi] these
tables for small values of $q$ (powers of $2$ or $3$) mainly by
a computer search on certain families. However, in many instances
the entries are large intervals of the form $[a,b]$ in which
$N_q(g)$ lies; here $b$ is the best upper bound known to him (from
Hasse-Weil-Serre or Ihara), while  $a$ indicates that he knows
that there exists a curve of genus $g$ with $a$ points over $\F_q$.
For completeness' sake we reproduce here his table. In order to improve
on his table one needs ways to construct  curves with many points.

The notion of generalized weight distribution in coding theory
applied to so-called trace codes suggested to us a construction of
curves with many points.  This led us to  considerable
improvements and extensions of the table of Wirtz. The
benefit  is here both ways: coding theory leads to curves with many
points and to curves with other interesting properties like
supersingular curves. On the other hand results in algebraic geometry
lead to the determination of generalized Hamming weights in coding
theory.

In this paper we concentrate on constructing curves with many
points. We will give the relation between curves and codes, explain
our method  and illustrate this with examples of the results. We will
also recall some other methods to find curves with many points. At the
end we give in addition to the table by Wirtz  two tables listing some
relatively good intervals or  actual values  for
$N_q(g)$.

It should  be clear to the reader that the story does not end here. As
the tables and the paper show there is ample room for improvement, both
theoretically and practically.

\bigskip
\noindent
\centerline {{\bf 1. The function } $N_q(g)$  }
\smallskip
\noindent
We begin with some notation. Let $q=p^m$ be a
power of a prime and let $\F_q$ be a finite field of cardinality $q$.
Let $C$ be an irreducible (complete non-singular) curve defined over
$\F_q$. We set $c(r)= \#C(\F_{q^r})$ and consider the zeta function
$$
Z(t)= {\rm exp}(\sum_{r=1}^{\infty} c(r) {t^r\over r}).
$$
According to Weil we can write
$$
Z(t)= {P_1(t) \over (1-t)(1-qt)},
$$
where $P_1(t) \in \Z[t]$ is of the form $P_1(t)= \prod_{i=1}^{2g}
(1-\alpha_it)$ with $\alpha_i \in \C$ algebraic integers of absolute
value
$\sqrt{q}$. Thus we have
$$
c(r)= q^r +1 - \sum_{i=1}^{2g} \alpha_i^r
$$
and this gives the famous Hasse-Weil bound
$$
| \# C(\F_q) - (q+1)| \leq 2g\sqrt{q},
$$
or more precisely
$$
| \# C(\F_q) - (q+1)| \leq [2g\sqrt{q}].
$$
(Here $[\cdot  ]$ denotes the greatest integer function.) Using
some algebraic number theory Serre improved this in 1983 to:
$$
| \# C(\F_q) - (q+1)| \leq g[2\sqrt{q}].
$$

We denote by $N_q(g)$ the maximum number of rational points on a curve
of genus $g$ over $\F_q$. Then by the above we have immediately
$$
N_q(g) \leq q+1+ g[2\sqrt{q}].
$$

For $g$ large with respect to $q$ there is a better bound.
The basic idea for this is due to Ihara [I] .  Suppose that
$\# C(\F_q)$ is large, then the $\alpha_i$ occurring above are near to
$-\sqrt{q}$ in $\C$. Hence  the $\alpha_i^2$ are near to $q$, so $\#
C(\F_{q^2}) $ is small, but we have $\# C(\F_q) \leq \# C(\F_{q^2})$.
Working this out with Cauchy-Schwarz gives a better bound. This
bound (the Ihara bound) is $N_q(g)\leq q+1 +[
(\sqrt{ (8q+1)g^2 + 4(q^2-q)g } -g)/2 ]$. Combined with Serre's bound
we find
$$
N_q(g) \leq \min \Bigl\{ q+1+g[2\sqrt{q} ], q+1 +[
(\sqrt{ (8q+1)g^2 + 4(q^2-q)g } -g)/2 ]\Bigr\}
$$
which is better than the Hasse-Weil bound for $g > (q-\sqrt{q})/2$.

Ihara's method was generalized by Drinfeld
and Vladut in [D-V] and by Serre. Serre uses  the idea of a `formule
explicite'. Take  a trigoniometric polynomial
$$
f= 1 + 2 \sum_{n\geq 1} u_n \cos n\theta
$$
which is even with real coefficients $u_n$ such  that  $f$
satisfies i)
$f(\theta) \geq 0$ for all $\theta \in \R$, ii) $u_n \geq 0$ for all
$n\geq 1$. We set
$$
\psi_1= \sum_{n\geq 1} u_nt^n.
$$
Then we have the estimate:
$$
N \leq a_f g + b_f,
$$
with
$$
a_f= {1 \over \psi_1(1/ \sqrt{q})} \qquad {\rm and } \qquad b_f = 1 +
{\psi_1(\sqrt{q}) \over \psi_1(1/ \sqrt{q})}.
$$
The game is to find the appropriate
function $f$, cf.\  [S 1]. Oesterl\'e
managed to find optimal choices for such $f$; we refer to [Sch]
for an exposition of  the Oesterl\'e algorithm.

Consider now the case where $g \leq (q-\sqrt{q})/2$. It was
conjectured by Stichtenoth and Xing that if $q$ is a square and  if $C$
is a {\sl maximal} curve of genus
$g$ over $\F_{q}$, i.e.  with
$\# C(\F_q)= q+1+2g\sqrt{q}$, then one should have
$$
g\leq (\sqrt{q}-1)^2/4 \quad {\hbox {\rm or } }\quad  g=
(q-\sqrt{q})/2.
$$
After considerable progress by Stichtenoth and Xing [S-X] this  was
recently proved  by Fuhrmann and  Torres [F-T]. This
excludes certain values for
$N_q(g)$ and improves the Hasse-Weil bound in some places (slightly).

Although the present paper deals with curves of small genus over a
small field it is natural to study the asymptotics of
$N_q(g)$. Define
$$
A(q)= \lim \sup N_q(g)/g \qquad {\rm as} \quad g \to \infty.
$$
The Drinfeld-Vladut bound gives
$$
A(q) \leq \sqrt{q}-1
$$
and then equality for $q=p^{2m}$  follows since in that case Ihara proves
$A(q) \geq \sqrt{q}-1$ using modular curves. Tsfasman, Vladut and Zink
used this to construct geometric Goppa codes which are better than the
Gilbert-Varshamov bound in coding theory, see [T-V-Z].  Recently, for
$q$ a square Garcia and Stichtenoth gave an  explicit construction with
repeated Artin-Schreier extensions of a sequence of curves for which
$N_q(g)/g$ goes to $\sqrt{q}-1$, cf.\ [G-S1]. 

Finally, we want 
to mention a recent survey paper [G-S2]  in this field of a different
nature.

\eject
\bigskip
\noindent
\centerline{\bf 2. Trace Codes, Higher Weights  and Curves  }
\bigskip
\noindent
We consider a (linear) code of length $n$:  an $\F_q$-linear subspace
$C$ of
$\F_q^n$. The elements of $C$ are called {\sl code words}.  A good deal
of information of  a code is stored in the {\sl weight distribution},  i.e.
in the polynomial
$$
A= \sum_{c \in C} X^{w(c)}= \sum_{i=0}^n A_i X^i  \quad (\in \Z[X])
$$
where $w(c)= \# \{i: c_i\neq 0\}$ is the  weight of the word
$c=(c_1,\ldots,c_n)$ and $A_i=\#\{ c\in C : w(c) = i\}$ is the
frequency of weight $i$.

Although this is a very simple situation, for some naturally defined
codes it turns out to be very difficult to determine the weight
distribution and determining the weight distribution is one of the
central problems of coding theory.

To convince the algebraic geometer of the difficulty of this, we
consider the classical {\sl Reed-Muller codes} $R(r,m)$  of order
$r$ over
$\F_p$. For the vector space
$$
P_r = \{ f \in \F_p[X_1, \ldots, X_m] : \deg(f) \leq r \}
$$
we have an evaluation map $\beta : P_r \to \F_p^n$ with $n=p^m$ given
by  $f \mapsto (f(v)_{v \in \F_p^m})$. The image $\beta(P_r)$ is
the code $R(r,m)$. The weight distribution of these codes is unknown
for $r\geq 3$. Indeed, e.g. for $r=3$  we need to know the distribution
of the number of points on all  cubic hypersurfaces (of a fixed space)
and this seems out of reach.
\smallskip

For reasons of simplicity and practice one is often first interested in
codes over a prime field, especially over $\F_2$ and $\F_3$.  In order
to get a code over the prime field from a code over an extension field
we have two methods:  applying the trace map or the restriction map. If
$\tilde C
\subset
\F_q^n$ is a (linear) code of length $n$ over $\F_q$, then the trace
code ${\rm Tr}(\tilde C )$ is obtained by applying the usual trace map
${\rm Tr}:\F_q\to
\F_p$ to all coordinates of all codewords in $\tilde C$. The trace
map is the most important for us.  Cyclic codes and many other
classical codes or their duals are trace codes.

A basic  example of trace codes that we shall exploit in the sequel
is obtained as follows.

Fix $h \geq 0$ and let
$$
{\cal R}_h = \{ R = \sum_{i=0}^h a_i X^{p^i}: a_i \in \F_q \}.
$$
This is a $\F_q$-vector space of {\sl additive} polynomials. We define
the code
$C_h$ as
$$
 C_h  =\{ c_R= ({\rm
Tr}(xR(x))_{x \in \F_q}): R \in {\cal R}_h \}.
$$
Since  additive polynomials behave as
$\F_p$-linear functions the expression
${\rm Tr}(xR(x))$ defines a quadratic form on the
$\F_p$-vector space $\F_q$. Thus the code $C_h$ is a subcode of the
classical binary Reed-Muller code $R(2,m)$.

Instead of the code $C_h$ above we can also consider the modified code
$C_h^*$ obtained from $C_h$ by puncturing at $x=0$, that is, deleting
the coordinate corresponding to $x=0$. It has the same weights as
$C_h$. For $h=1, 2$ the code $C_h^*$ is the dual of the classical $2$-
and $3$- error correcting BCH-codes ${\rm BCH}(2)$ and
${\rm BCH}(3)$.

Other examples of  classical codes which are trace codes and useful
for our purposes include  the dual Melas codes $M(q)^{\bot}$ of length
$q-1$, the words of  which are
$$
c_{a,b} = ({\rm Tr}(ax+b/x)_{x \in \F_q^*}) \quad {\rm for }\quad
a,b\in \F_q.
$$

\smallskip
A simple observation links  trace codes  with the theory of algebraic
curves.

The word $c_R= ({\rm Tr}(xR(x))_{x \in \F_q}) \in {\cal C}_h$ has
weight
$$
w(c_R)= \# \{ x \in \F_q : {\rm Tr}(xR(x))\neq 0\} = q - {1\over p}(\#
C_R(\F_q)-1),\eqno(1)
$$
with $C_R$ the (complete smooth) algebraic curve over $\F_q$ given by
$$
y^p-y=xR(x).
$$
The relation (1) is a simple consequence of the observation that
$$
{\rm Tr}(a)=0 \Longleftrightarrow {\hbox {\rm  there exists } } b \in
\F_q {\hbox {   \rm with } }  b^p-b = a.
$$
We find a correspondence $c_R \leftrightarrow C_R$ between non-zero
words of $C_h$ and algebraic curves. Thus codes lead to {\sl
families} of curves and the weight distribution is equivalent to the
distribution of the number of points on the curves in our family. It
explains why the weight distributions of such codes are hard to get
by. Note that words with low weight correspond to curves with many
points.

The dual Melas codes mentioned above lead to families of curves
$C_{a,b}$ given by
$$
y^p-y= ax+b/x
$$
with
$$
w(c_{a,b})= q-1 - {1 \over p} (\# C_{a,b}(\F_q)-2).\eqno(2)
$$
For $p=2$ this gives a
(universal) family of elliptic curves, while for $p=3$ we find a
family of genus $2$ curves which also leads to a universal family of
elliptic curves, see [G-V 1].  The codes
$C_h$ lead to {\sl supersingular} curves and these were studied in
[G-V 2].

\smallskip
We now come to the {\sl higher weight distribution} of a code. Although
the concept of the higher weights could have been introduced long ago,
interest for it arose only recently and was motivated by
applications, see [W]. Define for an $r$-dimensional  linear
subcode (subspace)
$D \subseteq C$ of a code $C$ of length $n$ over $\F_p$ the weight
$w(D)$ by
$$
w(D)= {1 \over p^r - p^{r-1}} \sum_{c \in D} w(c).
$$
An alternative definition is
$$
w(D)= n - \# \{ i : 1\leq i \leq n, c_i = 0
{\hbox { \rm for all } } c \in D\}.
$$

The $r$-th {\sl generalized Hamming weight} $d_r(C)$ of $C$ is
$$
d_r(C)= \min \{ w(D) : D \subseteq C, \dim (D) = r \}.
$$
So $d_1(C)$ is the usual minimum distance  of  $C$. The {\sl weight
hierarchy} of $C$ is defined as
$$
\{ d_r : 1\leq r \leq \dim(C) \}.
$$

For the trace codes above we saw that the  words in the code are
related to  curves in a family of curves. Similarly, the
subcodes are related to  fibre products
of such curves.  To be more precise, we start with a trace code $C$
with words of the form
$$
c_f= ({\rm Tr}(f(x))_{x \in \P^1(\F_q)-P}).
$$
Here $f$ runs over a finite-dimensional $\F_q$-subspace ${\cal L}$ of
the function field $\F_q(x)$ of $\P^1$. We require that the non-zero
elements of ${\cal L}$ are not of the form $g^p-g+a$ with $g \in
\F_q(x)$ and $a \in \F_q$. Furthermore, $P$ is a fixed subset of
$\P^1(\F_q)$ such that the elements of ${\cal L}- \{ 0 \}$ have poles
only in $P$. Note that for the codes $C_h$ we have $P= \{ \infty \}$
and for $M(q)^{\bot}$ the pole set $P= \{ 0 , \infty \}$.

Take an $r$-dimensional subcode $D$ of $C$ with basis $c_{f_1},
\ldots, c_{f_r}$. To these words correspond functions $f_1,
\ldots, f_r$ from ${\cal L}$ and they span an $r$-dimensional
$\F_p$-subspace ${\cal L}_D$ of ${\cal L}$. Let $C_{f_i}$ be the curve
defined by $y^p-y=f_i(x)$. Each
of these is an Artin-Schreier extension of $\P^1$ and thus admits a map
$\phi_i: C_{f_i} \to \P^1$. Then we can associate to the subcode $D$
the curve
$$
C_D = {\hbox {\rm Normalization of } } C_{f_1} \times_{\P^1}
\ldots \times_{\P^1} C_{f_r}.
$$
Up to isomorphism this does not depend on the chosen basis.

We then have the following
generalization of the relation between the weight of the subcode
$D$ and the number of
points of $C_D$, see [G-V 3]:
\smallskip
\noindent
\proclaim
(2.1) Proposition. The weight of the $r$-dimensional subcode $D$ of a
trace code
$C$ of length $n$  is related to the number of points of
$C_D$ by
$$
w(D)= n -  ( \# C_D(\F_q) - \sum_{Q \in  P}
p^{\epsilon(Q)})/p^r,\eqno(3)
$$
where $\epsilon(Q)= \dim_{\F_p} \{ f \in {\cal L}_D :
f {\hbox { \rm is regular at } } Q\}$.
\par

Of course, we need to know the trace of Frobenius for $C_D$. If
$\tau_f$ denotes the trace of Frobenius on $C_f$ (i.e. $\# C_f(\F_q) =
q+1- \tau_f$) then we have $$
(p-1) \tau_D = \sum_{f\in  {\cal L}_D-\{ 0 \} } \tau_f.\eqno(4)
$$
This can be proved by analyzing the number of points on $C_D$ in
relation to the number of points on the $C_f$, as we did in [G-V 6]. It
implies the existence of an isogeny $$
{\rm Jac}(C_D) \sim \prod_{f \in \P({\cal L}_D)} {\rm
Jac}(C_{f}),\eqno(5)
$$
where the summation is over a complete set of representatives for the
canonical action of
$\F_p^*$ on ${\cal L}_D- \{ 0 \}$. Note that we have isomorphisms
$C_f \cong C_{\lambda f}$ for any
$\lambda \in \F_p^*$, so that it makes sense to sum over $f$ in the
projective space of the $\F_p$-vector space ${\cal L}_D$.
Alternatively, using some Galois theory one can prove  (5) and then
deduce (4) from (5).

For  the genera of the curves we have
$$
(p-1) g(C_D)= \sum_{f \in {\cal L}_D - \{ 0 \}} g(C_f).\eqno(6)
$$

We see that if we  minimize the
weight of subcodes  $D$ of a fixed dimension we are maximizing the
number of $\F_q$-rational points on curves of the form $C_D$. Of
course, we need to control  the genera of the curves $C_f$ with
$f \in  {\cal L}_D - \{ 0 \} $.

\bigskip
\noindent
\centerline {\bf 3. Curves with Many Points }
\smallskip
\noindent
In this section we survey several methods to obtain linear spaces of
curves with many points.  These methods are closely related to coding
theory. By taking the fibre product corresponding to these linear
spaces we obtain new curves with many points which are distinguished
by the fact that they are given quite explicitly.  We illustrate each
method with some examples of the results. By combining the methods
we find results that go beyond those of [G-V 3-8]; there the reader
can find  a more detailed description of the methods and of some of
the results.
\bigskip
\noindent
{\bf Method I.} For $q=p^m$ with $m$ even the kernel of the canonical
$\F_p$-linear map
$$
\phi: {\cal R}_{m/2} \to C_{m/2}\quad{\hbox{\rm  defined by } } \quad
\phi (R) = ({\rm Tr}(xR(x))_{x \in \F_q})
$$
has dimension $m/2$ and consists of the additive polynomials
$R=ax^{\sqrt{q}}$ with $a \in \F_q^*$  satisfying $a^{\sqrt{q}} + a
=0$. For
$a\neq 0$ the correspoding curves, given by an affine equation
$$
y^p-y= ax^{\sqrt{q}+1},
$$
have genus $g=(p-1)\sqrt{q}/2$ and all have $pq+1$ points over $\F_q$.

Observe that these curves attain the Hasse-Weil upper bound.  Applying
the fibre product construction  and using the relations (4) and (6)
yields immediately the following result.
\smallskip
\noindent
\proclaim (3.1) Theorem. For $q=p^m$ with $m$ even and $1\leq r \leq
m/2
$ there is an $r$-dimensional subspace $L$ of curves over $\F_q$ with
$pq+1$ rational points. The corresponding fibre product yields a
maximal curve $C_L$ over $\F_q$ with
$$
g(C_L)= (p^r-1)\sqrt{q}/2 \quad { \hbox { \rm and } } \quad
\#C_L(\F_q)= p^rq+1.
$$
\par
\smallskip
For $q=p^m$ with $m$ odd the kernel of the canonical map
$$
\phi: {\cal R}_{(m+1)/2} \to C_{(m+1)/2}\quad{\hbox{\rm  defined by
}}\quad
\phi (R) = ({\rm Tr}(xR(x))_{x \in \F_q})
$$
has dimension $m$ and consists of the additive polynomials
$$
R= ax^{p^{(m+1)/2}}- (ax)^{p^{(m-1)/2}}\quad {\rm with }\quad a \in
\F_q.
$$
For $a\neq 0$ these curves have genus $(p-1)\sqrt{pq}/2$ and $pq+1$
points over $\F_q$. In the same way as above we find:
\smallskip
\noindent
\proclaim (3.2) Theorem.  For $q=p^m$ with $m$ odd and $1\leq r \leq m$
there is an $r$-dimensional subspace $L$ of curves over $\F_q$ with
$pq+1$ rational points. The corresponding fibre product yields a
 curve $C_L$ over $\F_q$ with
$$
g(C_L)= (p^r-1)\sqrt{pq}/2 \quad{ \hbox { \rm and } }\quad
\#C_L(\F_q)= p^rq+1.
$$
\par
\vskip 1.0pc
\noindent
This gives some immediate improvements of Wirtz's table:

\vs
\noindent
\medskip
\vbox{
\noindent{\bf (3.3) Example. $p=2$}
\medskip\centerline{\def\quad{\hskip 0.6em\relax}
\def\quod{\hskip 0.5em\relax }
\vbox{\offinterlineskip
\hrule
\halign{&\vrule#&\strut\quod\hfil#\quad\cr
height2pt&\omit&&\omit&&\omit&&\omit&&\omit&&\omit&\cr
&{\it Field} &&{\it dim}&&{\it genus}&&$\# C_r(\F_q)$&&{\it Wirtz}&&{\hbox{\it
upper bound}}&\cr
height2pt&\omit&&\omit&&\omit&&\omit&&\omit&&\omit&\cr
\noalign{\hrule}
height2pt&\omit&&\omit&&\omit&&\omit&&\omit&&\omit&\cr
\noalign{\hrule}
&$\F_8$&&$1$&&$2$&&$17$&&$18$&&&\cr
&\omit&&$2$&&$6$&&$33$&&$25-36$&&&\cr
&\omit&&$3$&&$14$&&$65$&&$65$&&&\cr
\noalign{\hrule}
&$\F_{32}$&&$1$&&$4$&&$65$&&$65-77$&&&\cr
&\omit&&$2$&&$12$&&$129$&&&&$165$&\cr
&\omit&&$3$&&$28$&&$257$&&$137-298$&&&\cr
&\omit&&$4$&&$60$&&$513$&&&&$542$&\cr
&\omit&&$5$&&$124$&&$1025$&&&&$1025$&\cr
\noalign{\hrule}
&$\F_{128}$&&$1$&&$8$&&$257$&&$257-305$&&&\cr
&\omit&&$2$&&$24$&&$513$&&&&$657$&\cr
&\omit&&$3$&&$56$&&$1025$&&&&$1361$&\cr
height2pt&\omit&&\omit&&\omit&&\omit&&\omit&&\omit&\cr }
\hrule}
}}
\medskip
\noindent

\medskip
\vbox{
\noindent{\bf (3.4) Example. $p=3$}
\medskip\centerline{\def\quad{\hskip 0.6em\relax}
\def\quod{\hskip 0.5em\relax }
\vbox{\offinterlineskip
\hrule
\halign{&\vrule#&\strut\quod\hfil#\quad\cr
height2pt&\omit&&\omit&&\omit&&\omit&&\omit&&\omit&\cr
&{\it Field} &&{\it dim}&&{\it genus}&&$\# C_r(\F_q)$&&{\it Wirtz}&&{\hbox{\it
upper bound}}&\cr
height2pt&\omit&&\omit&&\omit&&\omit&&\omit&&\omit&\cr
\noalign{\hrule}
height2pt&\omit&&\omit&&\omit&&\omit&&\omit&&\omit&\cr
\noalign{\hrule}
&$\F_3$&&$1$&&$3$&&$10$&&$10$&&&\cr
\noalign{\hrule}
&$\F_{27}$&&$1$&&$9$&&$82$&&$82-110$&&&\cr
&\omit&&$2$&&$36$&&$244$&&$184-319$&&&\cr
&\omit&&$3$&&$117$&&$730$&&&&$877$&\cr
height2pt&\omit&&\omit&&\omit&&\omit&&\omit&&\omit&\cr }
\hrule}
}}
\bigskip

Variations on this theme in which one  avails oneself of other
trace forms which vanish on $\F_q$   also yield quite a good harvest.
We mention a few fruitful ones:
\smallskip
\item{$\bullet$} ${\rm Tr}(x^t(ax^q-ax))$ for $a \in \F_{q=p^m}$
and $g.c.d.(t,p)=1$.
\smallskip
\item{$\bullet$}${\rm Tr}(ax^r-a^{p^t}x^s)$ for $a \in \F_{q=p^m}$
and $s\equiv rp^t (\bmod (p^m-1))$.

\smallskip
\item{$\bullet$}${\rm Tr}(x(ax^{p^{m-1}}+bx^{p^{m/2}}-a^px^p))$ for
$ a \in \F_{q=p^m}$ with $m\equiv 0 (\bmod 2)$ and $b\in
\F_{q}$ satisfying $b^{\sqrt{q}}+b=0$.

\smallskip
\noindent
{\bf (3.5) Example.} Using for $\F_4$ the building blocks ${\rm
Tr}(x^3)$,
${\rm Tr}(ax^3+ax)$ with $a \in \F_4-\F_2$, ${\rm Tr}(x^5+x^2)$,  and
${\rm Tr}(x^7+x^4)$ we find the following curves.
\medskip
\vbox{
\noindent{}
\medskip\centerline{\def\quad{\hskip 0.6em\relax}
\def\quod{\hskip 0.5em\relax }
\vbox{\offinterlineskip
\hrule
\halign{&\vrule#&\strut\quod\hfil#\quad\cr
height2pt&\omit&&\omit&&\omit&\cr
&{\it genus}&&$\# C(\F_4)$&&{\it upper bound }&\cr
height2pt&\omit&&\omit&&\omit&\cr
\noalign{\hrule}
height2pt&\omit&&\omit&&\omit&\cr
\noalign{\hrule}
&$5$&&$17$&&$18$&\cr
&$6$&&$17$&&$20$&\cr
&$8$&&$17$&&$24$&\cr
&$11$&&$25$&&$30$&\cr
&$13$&&$33$&&$33$&\cr
&$17$&&$33$&&$40$&\cr
&$18$&&$33$&&$42$&\cr
&$19$&&$33$&&$43$&\cr
&$27$&&$49$&&$56$&\cr
&$35$&&$49$&&$69$&\cr
&$37$&&$65$&&$72$&\cr
&$41$&&$65$&&$78$&\cr
&$42$&&$65$&&$80$&\cr
height2pt&\omit&&\omit&&\omit&\cr }
\hrule}
}}
\noindent
The upper bounds shown here are obtained with
Oesterl\'e's optimal choice algorithm.
In the sequel we call them Oesterl\'e bounds.

\bigskip
\noindent
{\bf Method II. } This method consists of studying for $R \in {\cal
R}_h$ (with $R\neq 0$) the quadratic form $Q = {\rm Tr}(xR(x))$ on the
$\F_p$-vector space $\F_q$ associated to words in $C_h$. The
bilinear form
$$
B(x,y)= {\rm Tr}(xR(y) + yR(x))
$$
is symmetric with radical
$$
W= \{ x \in \F_q : B(x,y) = 0 {\hbox { \rm for all } } y \in \F_q \}.
$$
We denote the dimension of the $\F_p$-space $W$  by $w$. For $p\neq 2$
the quadratic form $Q$ has rank $m-w$. For $p=2$ the form $B(x,y)$
is symplectic and this implies $w\equiv m (\bmod 2)$. Moreover, if
$a_h\neq 0$ then $0 \leq w \leq 2h$.

 From now on we shall restrict to the case $p=2$. In characteristic
$2$ we also consider the space
$$
W_0= \{ x \in W : Q(x) = 0 \} .
$$
We have $W=W_0$ or $\dim(W_0)= \dim (W)-1$. From the theory of
quadratic forms we deduce that the rank of $Q$ satisfies
$$
{\rm rank}(Q) =\cases{m-w & if $W=W_0$,\cr
m-w+1 & if $W\neq W_0$.\cr}
$$
If $W \neq W_0$  the number of points on $Q$ is $q/2$, while if
$W=W_0$ it is either $(q+\sqrt{2^wq})/2$ or $(q-\sqrt{2^wq})/2$. The
sign depends on the Witt-index  of the quadratic form.

To construct suitable quadrics  we start with elements
$a_i, b_i$ for
$i=1,\ldots, (m-w)/2$. We require that these $m-w$ elements are
$\F_2$-linearly independent  in
$\F_q$. Here $w$ is an integer with $w\equiv m (\bmod 2)$. Then we
claim that
$$
Q_{a,b}= \sum_{i=1}^{(m-w)/2} {\rm Tr}(a_ix){\rm Tr}(b_ix)
$$
is a quadratic form with $(q+\sqrt{q2^w})/2$ zeros in $\F_q$.
Note that ${\rm Tr}(ax)$ is a linear form so that the expression
defines indeed a quadratic form. By simple linear algebra one
transforms the quadratic form into the form $X_1X_2 + \ldots
X_{m-w-1}X_{m-w}$ with Witt-index $(m-w)/2$ and this implies the
result. So we can associate to
$$
\{ a_1,\ldots, a_{(m-w)/2}, b_1, \ldots, b_{(m-w)/2} \}
$$
a quadratic
form $Q_{(a,b)}$. We now try to arrange this so that $Q_{(a,b)}$ lies
in the code $C_h$ or rather $C_h^*$, see Section 2. That means that the
coefficients in
$$
\sum_{i=1}^{(m-w)/2} {\rm Tr}(a_ix){\rm Tr}(b_ix) =
\sum_{i=1}^{(m-w)/2} {\rm Tr}({\rm Tr}(a_ix)b_ix)
$$
of terms
$x^{{2^j}+1}$ with $j>h$ must disappear. The condition is:
\smallskip
\noindent
\proclaim (3.6) Proposition. If the elements $a_i, b_i \in \F_q$ with
$1
\leq i \leq (m-w)/2$ satisfy the system of equations
$$
\sum_{i=1}^{(m-w)/2} (a_i^{2^j}b_i + a_i b_i^{2^j}) =0 \eqno(7)
$$
for $j=h+1,\ldots , (m-1)/2$ then the word of length $q-1$ obtained
by evaluating $\sum_{i=1}^{(m-w)/2} {\rm Tr}(a_ix){\rm Tr}(b_ix)$ on
$\F_q^*$ is a codeword in $C_h^*$. \par
Note that these equations define so-called Deligne-Lusztig
varieties. For fixed $a_1,a_2,\ldots, a_{(m-w)/2}$ the words induced by
the solutions $b_1,\ldots,b_{(m-w)/2}$  of (7) form a linear space of
minimum weight words in $C_h^*$. In [G-V 8] we analyzed the dimension
of this space.  We find for example:

\smallskip
\proclaim  (3.7) Theorem. For $1 \leq r \leq m - 3$ with $m\geq 5$
and $m$ odd there exist curves
$C_r$ defined over $\F_{2^m}$ of genus $g (C_r)
= (2^r - 1) 2^{(m - 5) /2}$ and with $\# C_r (\F_{2^m}) = 2^m + 1 +
(2^r - 1) 2^{m - 2}$.
\par
\noindent
\smallskip
\noindent
{\bf  (3.8) Example.} For $m = 7$  Theorem (3.7) yields:
\smallskip

$\bullet$ {\sl There exist curves} $C_r$ {\sl
defined over} $\F_{128}$ {\sl with:}
\medskip
\medskip
\vbox{
\noindent{}
\medskip\centerline{\def\quad{\hskip 0.6em\relax}
\def\quod{\hskip 0.5em\relax }
\vbox{\offinterlineskip
\hrule
\halign{&\vrule#&\strut\quod\hfil#\quad\cr
height2pt&\omit&&\omit&&\omit&&\omit&\cr
&$r$&&$g(C_r)$ &&$\# C_r(\F_{128})$&&{\it Wirtz}&\cr
height2pt&\omit&&\omit&&\omit&&\omit&\cr
\noalign{\hrule}
height2pt&\omit&&\omit&&\omit&&\omit&\cr
\noalign{\hrule}
&$1$&&$2$&&$161$&&$172$&\cr
&$2$&&$6$&&$225$&&$225-261$&\cr
&$3$&&$14$&&$353$&&$289-437$&\cr
&$4$&&$30$&&$609$&&$369-789$&\cr
height2pt&\omit&&\omit&&\omit&&\omit&\cr }
\hrule}
}}
\medskip
\noindent

\medskip
\proclaim (3.9) Theorem. For $1 \leq r \leq (m-2)/2$ with $m$ even
and
$m\geq 4$ there exist curves $C_r$ defined over $\F_{2^m}$ of genus
$g(C_r)=(2^r-1)2^{(m-4)/2}$ and with $\# C_r(\F_{2^m})= 2^m+1 +
(2^r-1)2^{m-1}$; these curves attain the Hasse-Weil bound.
\par
\smallskip

\smallskip
\noindent
{\bf  (3.10) Example.} For $m = 6$ (resp. $m=8$) this Theorem
yields:

$\bullet$ {\sl There exist curves} $C_r$ {\sl
defined over} $\F_{64}$  (resp. over $\F_{256}$) {\sl as indicated:}
\medskip

\smallskip
\medskip
\vbox{
\noindent{}
\medskip\centerline{\def\quad{\hskip 0.6em\relax}
\def\quod{\hskip 0.5em\relax }
\vbox{\offinterlineskip
\hrule
\halign{&\vrule#&\strut\quod\hfil#\quad\cr
height2pt&\omit&&\omit&&\omit&&\omit&&\omit&&\omit&\cr
&{\it Field} && $r$&&$g(C_r)$ &&$\# C_r(\F_{q})$&&{\it Wirtz}&&{\it upper
bound}&\cr
height2pt&\omit&&\omit&&\omit&&\omit&&\omit&&\omit&\cr
\noalign{\hrule}
height2pt&\omit&&\omit&&\omit&&\omit&&\omit&&\omit&\cr
\noalign{\hrule}
&$\F_{64}$&&$1$&&$2$&&$97$&&$97$&&&\cr
&\omit&&$2$&&$6$&&$161$&&$155-161$&&&\cr
\noalign{\hrule}
&$\F_{256}$&&$1$&&$4$&&$385$&&&&$385$&\cr
&\omit&&$2$&&$12$&&$641$&&&&$641$&\cr
&\omit&&$3$&&$28$&&$1153$&&&&$1153$&\cr
height2pt&\omit&&\omit&&\omit&&\omit&&\omit&&\omit&\cr }
\hrule}
}}
\medskip
\noindent

\bigskip
\noindent
{\bf Method III. } The third approach also arose from quadratic
forms associated to the codes
$C_h$. In studying the weight distribution of $C_h$  in [G-V 2]  it was
convenient first to determine for fixed $R\in {\cal R}_h$ with $R \neq
0$ the weight distribution in the $1$-dimensional family ${\cal F}_R$
of words
$$
({\rm Tr}(xR(x)+bx))_{x \in \F_q}\quad {\rm with} \quad b \in \F_q.
$$
The weight distribution of ${\cal F}_R$ follows from the theory of
quadratic forms over finite fields. Only few weights occur and the
values of $b$ which provide minimum weight words or curves with many
points happen to lie in quadratic spaces. Using the Witt-index, which
determines the dimension of a maximal totally singular subspace in the
quadratic space, we find linear or affine spaces of curves with many
points. Again we have a situation in which we can profitably apply
the fibre product construction. From [G-V 7] we recall a result for
fields
$\F_{q=p^m}$ with $p$ an odd prime.
\medskip
\noindent
\proclaim (3.11) Theorem. i) If $m$ is odd there exists for $1 \leq r
\leq (m-1)/2$ a curve $C_r$ of genus $g=p^r(p-1)/2$ with $\# C_r(\F_q)=
q+1 + p^r \sqrt{pq}$.
ii) If $m$ is even there exists for $1 \leq r \leq m/2$ a curve $C_r$
of genus $g=p^r(p-1)/2$ which attains the Hasse-Weil upper bound: $\#
C_r(\F_q)= q+1+p^r(p-1)\sqrt{q}$
\par
\smallskip
If we specialize to characteristic $3$ and use some specific
properties of the occurring situations we derive the following
results.
\smallskip
\noindent
\proclaim (3.12) Theorem. i) For $q=3^m$ with $m$ odd, $m \geq 3$
there exists a curve of genus $4$ over $\F_q$ with $q+1+4\sqrt{3q}$
points over $\F_q$.
ii) For $q=3^m$ with $m\equiv 0( \bmod 4)$ there exists a curve of
genus $12$ over $\F_q$ which realizes the Hasse-Weil upper bound.
\par
\smallskip
\noindent
{\bf (3.13) Example.}
\medskip
\vbox{
\noindent{}
\medskip\centerline{\def\quad{\hskip 0.6em\relax}
\def\quod{\hskip 0.5em\relax }
\vbox{\offinterlineskip
\hrule
\halign{&\vrule#&\strut\quod\hfil#\quad\cr
height2pt&\omit&&\omit&&\omit&&\omit&&\omit&\cr
&{\it Field} &&{\it genus} &&$\# C_r(\F_{q})$ &&{\it Wirtz}&&{\it upper
bound}&\cr
height2pt&\omit&&\omit&&\omit&&\omit&&\omit&\cr
\noalign{\hrule}
height2pt&\omit&&\omit&&\omit&&\omit&&\omit&\cr
\noalign{\hrule}
&$\F_{27}$&&$3$&&$55$&&$55-58$&&&\cr
&\omit&&$4$&&$64$&&$55-68$&&&\cr
\noalign{\hrule}
&$\F_{81}$&&$12$&&$298$&&$226-298$&&&\cr
\noalign{\hrule}
&$\F_{243}$&&$3$&&$325$&&&&$337$&\cr
&\omit&&$4$&&$352$&&&&$368$&\cr
&\omit&&$9$&&$487$&&&&$523$&\cr
height2pt&\omit&&\omit&&\omit&&\omit&&\omit&\cr }
\hrule}
}}
\medskip
\noindent
\bigskip
\noindent
{\bf Method IV.} This method uses results from coding theory,
especially results about weight hierarchies of codes. From formula (3)
which describes the relation between the weight of a subcode $D$ and
$\# C_D(\F_q)$ we derive immediately:
\smallskip
\noindent
\proclaim (3.14) Proposition. i) For an $r$-dimensional subcode $D$ of the
punctured code $C_h^*$ of weight $w(D)$ the corresponding fibre
product curve $C_D$ satisfies:
$$
g(C_D)= {1 \over p-1} \sum_{f \in L_D - \{ 0 \} } g(C_f) \leq
(p^r-1)p^h/2
$$
and
$$
\# C_D(\F_q) = p^r(q-w(D))+1.
$$
ii) For an $r$-dimensional subcode $D$ of the dual Melas code
$M(q)^{\bot}$ of weight $w(D)$ the corresponding curve $C_D$
satisfies
$$
g(C_D)= {1 \over p-1} \sum_{f \in L_D - \{ 0 \} } g(C_f) \leq p^r-1
$$
and
$$
\# C_D(\F_q)  \geq p^r(q-1-w(D))+2.\eqno(8)
$$
Equality holds in (8) if $D$ contains no code words corresponding to
rational curves.
\par
\smallskip
An interesting situation appears if $D$ is a subcode of which the
weight equals the $r$-th generalized Hamming weight of the code. To
determine the value of the genus $g(C_D)$ one has to investigate
carefully the words occurring in $D$. In the following examples we
apply Proposition (3.14).
\smallskip
\noindent
{\bf (3.15) Example.} For $p=2$ and $q=2^m$ with $m\geq 5$ we proved
in [G-V 5] that the second generalized Hamming weight of
$BCH(2)^{\bot}$ of length $q-1$ satisfies
$$
d_2(BCH(2)^{\bot}) = d_2(C_1^*)= {3\over 2} d_1(BCH(2)^{\bot})=
\cases{ 3(q-\sqrt{2q})/4 & for $m$ odd,\cr
3(q-2\sqrt{q})/4 & for $m$ even.\cr}
$$
The curves corresponding to the $2$-dimensional subcode in which the
non-zero words all have minimum weight, have genus $1$. We deduce
from  Proposition (3.14):
\smallskip
\noindent
\proclaim (3.16) Theorem. For $m\geq 5$ there exists a curve $C$
over
$\F_{q=2^m}$ of genus $3$ with $\# C(\F_q)= q+1+3\sqrt{2q}$ for $m$
odd, while $C$ attains the Hasse-Weil upper bound for $m$ even.
\par
\smallskip
\noindent
{\bf (3.17) Example.} Since the code $BCH(2)$ of length $7$ is the
very simple repetition code $\F_2 \cdot {\underline{1}} \in \F_2^7$ the
weight hierarchy of $BCH(2)^{\bot}$ of length $7$ is
$$
\{ 2, 3, 4, 5, 6, 7\}.
$$
 From  Proposition (3.14) and an analysis of the genera of the curves
belonging to the subcodes we get

\medskip
\vbox{
\noindent{}
\medskip\centerline{\def\quad{\hskip 0.6em\relax}
\def\quod{\hskip 0.5em\relax }
\vbox{\offinterlineskip
\hrule
\halign{&\vrule#&\strut\quod\hfil#\quad\cr
height2pt&\omit&&\omit&&\omit&&\omit&&\omit&&\omit&\cr
&{\it Field} && $r$ && $d_r(BCH(2)^{\bot})$&&{\it genus } &&$\#
C_r(\F_q)$&&{\it Wirtz}&\cr
height2pt&\omit&&\omit&&\omit&&\omit&&\omit&&\omit&\cr
\noalign{\hrule}
height2pt&\omit&&\omit&&\omit&&\omit&&\omit&&\omit&\cr
\noalign{\hrule}
&$\F_8$ &&$1$&&$2$&&$1$&&$13$&&$14$&\cr
&\omit&&$2$&&$3$&&$3$&&$21$&&$24$&\cr
&\omit&&$3$&&$4$&&$6$&&$33$&&$25-36$&\cr
&\omit&&$3$&&$4$&&$7$&&$33$&&$25-39$&\cr
height2pt&\omit&&\omit&&\omit&&\omit&&\omit&&\omit&\cr }
\hrule}
}}
\medskip
\noindent

\smallskip
\noindent
{\bf (3.18) Example.} For $p=2$ and $h=[m/2]$ one derives by equating
dimensions that the code $\F_2 \cdot {\underline 1} + C_h$ is the
binary second order Reed-Muller code $R(2,m)$ of length $2^m$.
Elementary properties of the minimum weight words in $R(2,m)$ imply
that
$$
d_r= (2^r-1)d_1= (2^r-1) 2^{m-2} / 2^{r-1} \quad {\rm for }
\quad 1\leq r < m.
$$
For odd $m$ the curves associated to minimum weight words have genus
$2^{h-1}$. If we apply Proposition (3.14) to an $r$-dimensional
subcode of minimum weight words we find:
\smallskip
\noindent
\proclaim (3.19) Theorem. For odd $m$ with $m\geq 3$ and $r$ with
$1\leq r < m$ there exists a curve $C_r$ over $\F_{q=2^m}$ of genus $g=
  (2^r-1)\sqrt{q/8}$ with $\# C_r(\F_q) = q+1 +
(2^r-1)q/2$.
\par
\medskip
Theorem  (3.19) yields for instance:

\medskip
\vbox{
\noindent{}
\medskip\centerline{\def\quad{\hskip 0.6em\relax}
\def\quod{\hskip 0.5em\relax }
\vbox{\offinterlineskip
\hrule
\halign{&\vrule#&\strut\quod\hfil#\quad\cr
height2pt&\omit&&\omit&&\omit&&\omit&&\omit&\cr
&{\it Field} && $r$ &&{\it genus } &&$\# C_r(\F_q)$&&{\it Wirtz}&\cr
height2pt&\omit&&\omit&&\omit&&\omit&&\omit&\cr
\noalign{\hrule}
height2pt&\omit&&\omit&&\omit&&\omit&&\omit&\cr
\noalign{\hrule}
&$\F_{32}$&&$1$&&$2$&&$49$&&$53$&\cr
&\omit&&$2$&&$6$&&$81$&&$73-99$&\cr
&\omit&&$3$&&$14$&&$145$&&$113-187$&\cr
&\omit&&$4$&&$30$&&$273$&&$129-313$&\cr
height2pt&\omit&&\omit&&\omit&&\omit&&\omit&\cr }
\hrule}
}}
\medskip
\noindent

\medskip
We conclude with an example in which Melas codes and dual Melas
codes are central (i.e. involving non-supersingular elliptic curves).
\smallskip
\noindent
{\bf (3.20) Example.} From [G-V 4] one can derive the weight hierarchy
of the dual Melas code $M(16)^{\bot}$ of length $15$:
$$
\{ 4, 6, 8 , 9, 11, 12, 14, 15 \}.
$$
For $1\leq r \leq 4$ one can find an $r$-dimensional subcode
realizing the first four generalized Hamming weights in which all the
non-zero words belong to elliptic curves. Application of Proposition
(3.14)  provides us with the following table.
\smallskip

\medskip
\vbox{
\noindent{}
\medskip\centerline{\def\quad{\hskip 0.6em\relax}
\def\quod{\hskip 0.5em\relax }
\vbox{\offinterlineskip
\hrule
\halign{&\vrule#&\strut\quod\hfil#\quad\cr
height2pt&\omit&&\omit&&\omit&&\omit&&\omit&&\omit&\cr
&{\it Field} && $r$ && $d_r(M(16)^{\bot})$&&{\it genus } &&$\#
C_r(\F_q)$&&{\it Wirtz}&\cr
height2pt&\omit&&\omit&&\omit&&\omit&&\omit&&\omit&\cr
\noalign{\hrule}
height2pt&\omit&&\omit&&\omit&&\omit&&\omit&&\omit&\cr
\noalign{\hrule}
&$\F_{16}$&&$1$&&$4$&&$1$&&$24$&&$25$&\cr
&\omit&&$2$&&$6$&&$3$&&$38$&&$38$&\cr
&\omit&&$3$&&$8$&&$7$&&$58$&&$49-70$&\cr
&\omit&&$4$&&$9$&&$15$&&$98$&&$49-113$&\cr
height2pt&\omit&&\omit&&\omit&&\omit&&\omit&&\omit&\cr }
\hrule}
}}
\medskip
\noindent

\smallskip
\noindent
For other results where we use Melas codes and dual Melas codes we
refer to [G-V 3] and [G-V 7].

\eject

We finish  this section by observing that combining  the methods,
especially I and IV, can be rewarding. We mention some results for
$\F_8$.
\medskip
\vbox{
\noindent{\bf (3.21) Example.}
\medskip\centerline{\def\quad{\hskip 0.6em\relax}
\def\quod{\hskip 0.5em\relax }
\vbox{\offinterlineskip
\hrule
\halign{&\vrule#&\strut\quod\hfil#\quad\cr
height2pt&\omit&&\omit&&\omit&\cr
&{\it genus}&&$\# C(\F_8)$&&{\it upper bound }&\cr
height2pt&\omit&&\omit&&\omit&\cr
\noalign{\hrule}
height2pt&\omit&&\omit&&\omit&\cr
\noalign{\hrule}
&$5$&&$29$&&$32$&\cr
&$11$&&$41$&&$54$&\cr
&$13$&&$49$&&$61$&\cr
&$18$&&$57$&&$77$&\cr
&$22$&&$65$&&$89$&\cr
&$23$&&$65$&&$92$&\cr
&$27$&&$81$&&$103$&\cr
&$29$&&$97$&&$109$&\cr
&$38$&&$113$&&$135$&\cr
height2pt&\omit&&\omit&&\omit&\cr }
\hrule}
}}
\medskip
\noindent
 The upper bound in the table is Oesterl\'e's.
\bigskip
\noindent
\centerline {\bf 4.  Final Remarks }
\smallskip
\noindent
There are of course other ways of producing curves with many points.
Besides a direct computer search on certain families as performed
by Wirtz one could use specific curves like modular curves. In our
tables we have incorporated some values obtained from modular curves.
But more in line with the spirit of this paper,  one can also  use
coverings of curves of genus
$ > 0$ with many points to produce new curves with many points. Note
that we used so far only Artin-Schreier covers of $\P^1$, but there is
no reason to stick to this restriction.  For example, let $C$ be a
complete non-singular curve over $\F_q$ with genus $g$ whose
$\F_q$-rational points are $\{ P_1, \ldots ,P_m\}$. Furthermore, let
$Q=Q_n$ be a divisor of degree $n$ defined over $\F_q$ whose support
contains none of the $P_i$. Let $L \subset \F_q(C)$ be an
$r$-dimensional  $\F_p$-vector space of functions contained in
$$
\{ f
\in \F_q(C)^* : {\rm div}(f)
\geq -Q \} \cup  \{  0  \}.
$$
We shall assume that $L-\{ 0 \}$ does not contain elements of the form
$g^p-g+a$ with $g \in \F_q(C)$ and $a \in \F_q$. For $f \in L$ we can
form the covering
$$
C_f \to C, \quad {\hbox {\rm  given by  } } \quad y^p-y =f.
$$
This is a $p$-fold cover of $C$ ramified only in the poles of $f$. As
before we can also define $C_L$ using  the fibre
product $C_{f_1} \times_C \ldots \times_C C_{f_r}$ for a basis
$f_1,\ldots, f_r$ of $L$. The number of points can be expressed as
$$
\# C_L(\F_q) = p^r \cdot \# \{ P_i : {\rm Tr}(f(P_i))=0 {\hbox {
\rm  for all } } f
\in L \}.
$$
The genus and trace of Frobenius can be computed from the isogeny
$$
{\rm Jac}(C_L) \sim {\rm Jac}(C) \times \prod_{f \in \P(L)}
P(C_f/C),
$$
where the product is over a complete set of representatives in
$L-\{ 0 \}$ for the natural action of $\F_p^*$ on $L-\{ 0 \}$.
Moreover,  $P(C_f/C)$ stands for the Prym variety of $C_f \to C$,
which is the  connected component of the kernel of the norm map ${\rm
Nm}: {\rm Jac}(C_f) \to {\rm Jac}(C)$. We find
$$
\tau_{C_L} = \tau_C + \sum_{f \in \P(L)} (\tau_{C_f} - \tau_C),
$$
a relation that can also be obtained by counting points on the
curves $C_L$ and $C_f$ as in [G-V 6], and
$$
g(C_L)= g(C) + \sum_{f \in \P(L)}(g(C_f) - g(C)).
$$
We need the genera of the $C_f$. If we assume that $f$ has $n_f$
simple poles we find
$$
g(C_f) = g(C) + (p-1)(g(C) + n_f -1).
$$

In order to construct curves with many points one starts with a curve
$C$ with many rational points over $\F_q$ and one considers $C_L$
such that ${\rm Tr}(f(P))=0$ for all (or at least for many) rational
points of $C$. So we are considering the trace codes associated to
geometric Goppa codes.

The simplest case is of course where
$q=p$ so that
${\rm Tr}(f(P))=0$ simplifies to $f(P)=0$.
\smallskip
\noindent
{\bf (4.1) Example.} Let $E$ be an elliptic curve over $\F_2$ with $5$
rational points $P_1, \ldots, P_5$. Let $Q_n$ be a divisor of degree
$n\geq 5$ over $\F_2$. Then $\dim_{\F_2}(L(Q_n-(P_1+\ldots +P_5)))=
n-5$ for $n \geq 6$ and $0$ or $1$  for $n=5$. We find for the choice
$L=L(Q_n-(P_1+\ldots +P_5))$ the formulas
$$
g(C_L) = 1 + (2^{n-5}-1)n, \qquad \# C_L(\F_2)= 5 \cdot 2^{n-5}\qquad
{\rm for}\quad n\geq 6.
$$
If we take for $n=5$ a divisor $Q_5$ which is linearly equivalent to
$\sum_{i=1}^5 P_i$ then we get a curve $C_L$ with $g=7$ and  $\#
C(\F_2)=10$.

One should compare this with the example that Serre gives in [S 1]
using class field theory. He defines a series of curves $C_n$ over
$\F_2$ with
$$
g(C_n)= 1 + (2^{n-4}-1)n, \qquad \# C_n(\F_2)= 5 \cdot 2^{n-4}\qquad
{\rm for}\quad n\geq 5.
$$
and thus finds a curve of genus $6$ with $10$ points. The reason for
the non-optimality of our example comes from the fact that we require
all functions $f\in L$ to have simple poles in $Q_n$, which excludes
using functions that have multiple poles there but which are
Artin-Schreier equivalent to a function $\tilde f$ with simple poles at
$Q_n$.

This phenomenon also explains why  often with the methods of
Section 3 we find curves with a large number of points, but not the
maximum possible. On the other hand, our curves are given 
explicitly.

\smallskip
\noindent
{\bf (4.2) Example.} Variations of Serre's method provide other 
good curves. Take a curve of genus $2$ over $\F_2$ with
$6$ points.  By taking a 
$4$-fold cover ramified  in a point of degree $7$   we find a curve
$C_{26}$ of genus
$g=26$ with $\# C_{26}(\F_2)=24$.

In this context we  wish to mention a result which was communicated to us by
Ren\'e Schoof. It follows from class field theory and yields good
curves over $\F_2$.
\smallskip
\noindent
\proclaim (4.3) Theorem. Let $C$ be a curve of genus $g$ over $\F_2$
which has $n$ rational points and a point of degree $m$ with $m\geq
n$. For $1 \leq l \leq m-n+1$ there exists a fibre product $\tilde
C$ over $C$ of degree $2^l$ with
$$
\tilde g \leq (2^l-1)(g+m-1) +g
$$
and with
$$
\# \tilde C(\F_2)= 2^ln.
$$
If the inequality for the
genus is strict then there exists an Artin-Schreier cover $C_f$ of $C$
with $g(C_f) = 2g-1$ and
$\#C_f(\F_2) = 2n$.
\par
\smallskip
Finally, Drinfeld modules provide a theory of explicit class fields of 
function fields.  Very recently Niederreiter and Xing obtained along 
these lines some good curves over $\F_3$ and $\F_4$ which we include here in
our tables, cf. [N-X] and [X].

\bigskip
\noindent
\centerline{\bf Tables}
\smallskip
\noindent
We reproduce here the table by Wirtz (in which we corrected a few
misprints) and give two other tables.  In his table Wirtz mentions not
only  the intervals for
$N_q(g)$, but also indicates for many values where they can be found
in the literature or from which curves or families they are obtained.
We refrain from giving this information. In many places Wirtz mentions
instead of an interval $a-b$ only an upper bound $-b$, for which he
takes  the Ihara-Serre bound. Note that often the Oesterl\'e bound
is better.

In the other two tables, one for characteristic $2$ and one for
characteristic $3$, we give for $g \leq 50$ an interval $[a,b]$ in
which
$N_q(g)$ lies. For $b$ we take either the Ihara-Serre bound (taking
into account [F-T]) or the Oesterl\'e bound if that one is better. For
$g\leq 50$ these bounds coincide if $q\geq 27$.

We only enter new values in the table if we consider them
reasonable. We have taken this to mean that
$a\geq [b/\sqrt{2}]$. This ratio is inspired by the fact that the
Ihara bound implies
$$
A(q)\leq  {1 \over 2}(\sqrt{(8q+1)}-1) < \sqrt{2q}
$$
while we know that $A(q)\leq \sqrt{q}-1< \sqrt{q}$. We would like to
stress here that in no way we aimed at completeness and just entered
values or intervals that we happen to know. A systematic search using
our methods should lead to  improvements and extensions of the
tables.
\vfill\eject
%
%
%
%
%
%
%
\def\quad{\hskip 0.3em\relax}
\def\quod{\hskip 0.3em\relax }
\font\tablefont=cmr8
\def\vhop{height2pt&\omit&&\omit&&\omit&&\omit&&\omit&&\omit&&\omit&&
      \omit&&\omit&&\omit&&\omit&&\omit&\cr}
\noindent{\bf The Table of Wirtz.}
\medskip
\line{\hss\vbox{
\tablefont
\lineskip=1pt
\baselineskip=10pt
\lineskiplimit=0pt
\setbox\strutbox=\hbox{\vrule height .7\baselineskip
                                depth .3\baselineskip width0pt}%
\offinterlineskip
\hrule
\halign{&\vrule#&\strut\quod\hfil#\quad\cr
\vhop
&$g\backslash q$&&2&&3&&4&&8&&9&&16&&27&&32&&64&&81&&128&\cr
\vhop
\noalign{\hrule}
\vhop
&1&&5&&7&&9&&14&&16&&25&&38&&44&&81&&100&&150&\cr
&2&&6&&8&&10&&18&&20&&33&&48&&53&&97&&118&&172&\cr
&3&&7&&10&&14&&24&&28&&38&&55--58&&62--66&&113&&136&&184--195&\cr
&4&&8&&12&&15&&--29&&26--30&&41--49&&55--68&&65--77&&129&&154&&200--217&\cr
&5&&9&&--15&&--18&&--32&&--36&&--57&&55--78&&--88&&101--145&&152--172&&202--239&\cr
&6&&10&&--17&&17--21&&25--36&&29--40&&65&&64--88&&73--99&&137--161&&190&&225--261&\cr
&7&&10&&--19&&--23&&25--39&&28--43&&49--70&&64--98&&81--110&&177&&160--208&&241--283&\cr
&8&&11&&--21&&--25&&--43&&--47&&--76&&--108&&--121&&--193&&--226&&257--305&\cr
&9&&12&&--23&&--28&&33--47&&37--51&&49--81&&82--118&&81--132&&173--209&&244&&209--327&\cr
&10&&12--13&&--25&&--30&&--50&&--55&&61--87&&--128&&--143&&139--225&&226--262&&--349&\cr
\noalign{\hrule}
&11&&13--14&&&&&&&&&&&&&&&&&&&&&\cr
&12&&14--15&&&&&&37--57&&37--63&&61--97&&91--148&&&&257&&226--298&&193--393&\cr
&13&&14--15&&&&&&&&&&&&&&&&169--273&&163--316&&&\cr
&14&&15--16&&&&&&65&&43--70&&81--108&&&&113--187&&193--289&&218--334&&289--437&\cr
&15&&17&&&&&&43--68&&&&49--113&&&&97--196&&197--305&&&&369--459&\cr
&16&&16--18&&&&&&&&46--78&&&&100--178&&&&199--321&&370&&&\cr
&17&&17--18&&&&&&&&&&&&&&&&&&&&&\cr
&18&&18--19&&&&&&&&&&&&107--192&&&&&&&&201--525&\cr
&19&&20&&&&&&&&&&&&&&&&&&&&&\cr
&20&&19--21&&&&&&&&&&&&&&&&&&&&&\cr
\noalign{\hrule}
&21&&21&&&&&&57--89&&64--97&&89--145&&&&121--244&&257--401&&250--460&&297--591&\cr
&24&&&&&&&&&&64--108&&&&136--235&&&&242--449&&298--514&&&\cr
&26&&&&&&&&&&&&&&&&157--283&&&&&&&\cr
&27&&&&&&&&&&&&&&&&&&&&313--568&&&\cr
&28&&&&&&&&65--114&&73--123&&97--181&&136--263&&137--298&&513&&370--586&&409--745&\cr
&30&&&&&&&&&&&&81--192&&157--277&&129--313&&257--536&&&&369--789&\cr
&35&&&&&&&&&&&&105--218&&&&&&&&&&&\cr
&36&&&&&&&&&&82--154&&&&184--319&&&&&&730&&&\cr
&39&&33&&&&&&&&&&&&&&&&&&&&&\cr
&42&&&&&&&&&&&&113--254&&157--360&&&&&&&&&\cr
&45&&&&&&&&&&&&&&&&193--428&&401--706&&&&561--1119&\cr
&48&&&&&&&&&&&&&&244--402&&&&&&&&&\cr
&49&&&&&&&&&&&&129--291&&&&&&&&&&&\cr
&50&&40&&&&&&&&&&&&&&&&&&&&&\cr
\vhop
}
\hrule
}
\hss}
\vfill\eject
%
%
%
%
%
%
%
%
\def\quad{\hskip 0.6em\relax}
\def\quod{\hskip 0.6em\relax}
\def\vhop{
    height2pt&\omit&&\omit&&\omit&&\omit&&\omit&&\omit&&\omit&&\omit&\cr}
\noindent{\bf Table p=2.}
$$
\vcenter{
\tablefont
\lineskip=1pt
\baselineskip=10pt
\lineskiplimit=0pt
\setbox\strutbox=\hbox{\vrule height .7\baselineskip
                                depth .3\baselineskip width0pt}%
\offinterlineskip
\hrule
\halign{&\vrule#&\strut\quod\hfil#\quad\cr
\vhop
&$g\backslash q$&&2&&4&&8&&16&&32&&64&&128&\cr
\vhop
\noalign{\hrule}
\vhop
&1&&5&&9&&14&&25&&44&&81&&150&\cr
&2&&6&&10&&18&&33&&53&&97&&172&\cr
&3&&7&&14&&24&&38&&63--66&&113&&190--195&\cr
&4&&8&&15&&25--29&&45--49&&70--77&&129&&200--217&\cr
&5&&9&&17--18&&29--32&&49--56&&73--88&&130--145&&227--239&\cr
&6&&10&&20&&33--36&&65&&81--99&&161&&225--261&\cr
&7&&10&&21--22&&33--39&&58--70&&89--110&&177&&241--283&\cr
&8&&11&&21--24&&31--43&&61--76&&&&&&257--305&\cr
&9&&12&&26&&36--47&&62--81&&&&173--209&&241--327&\cr
&10&&13&&27--28&&&&65--87&&&&193--225&&289--349&\cr
\noalign{\hrule}
&11&&14&&25--30&&41--54&&75--92&&&&201--241&&&\cr
&12&&14--15&&28--31&&47--57&&&&129--165&&257&&321--393&\cr
&13&&15&&33&&49--61&&97--103&&&&199--273&&&\cr
&14&&15--16&&29--35&&65&&97--108&&145--187&&&&353--437&\cr
&15&&17&&33--37&&&&98--113&&&&257--304&&369--459&\cr
&16&&16--18&&36--38&&56--71&&93--118&&&&&&&\cr
&17&&17--18&&40&&61--74&&&&&&&&&\cr
&18&&18--19&&33--42&&57--77&&99--129&&&&&&&\cr
&19&&20&&36--43&&&&&&&&&&&\cr
&20&&19--21&&&&61--83&&121--140&&&&&&&\cr
\noalign{\hrule}
&21&&21&&40--47&&&&125--145&&&&&&&\cr
&22&&21--22&&33--48&&66--89&&129--150&&&&&&&\cr
&23&&22--23&&36--50&&68--92&&123--155&&&&&&&\cr
&24&&20--23&&&&&&&&&&&&513--657&\cr
&25&&24&&&&&&&&&&&&&\cr
&26&&24--25&&55&&&&129--171&&&&&&&\cr
&27&&22--25&&49--56&&81--103&&&&&&401--497&&&\cr
&28&&24--26&&&&91--106&&129--181&&257--298&&513&&577--745&\cr
&29&&25--27&&&&97--109&&&&&&&&&\cr
&30&&24--27&&&&&&&&273--313&&&&609--789&\cr
\noalign{\hrule}
&31&&27--28&&&&&&&&&&&&&\cr
&32&&26--29&&&&&&&&&&&&&\cr
&33&&28--29&&65--66&&&&&&&&&&&\cr
&34&&27--30&&57--68&&&&&&&&&&&\cr
&35&&28--31&&49--69&&&&&&253--352&&&&&\cr
&36&&30--31&&&&105--130&&&&&&&&&\cr
&37&&28--32&&65--72&&121--132&&&&&&&&&\cr
&38&&28--33&&&&113--135&&&&&&&&&\cr
&39&&33&&&&117-138&&&&&&&&&\cr
&40&&32--34&&&&&&&&&&&&&\cr
\noalign{\hrule}
&41&&30--35&&65--78&&&&&&&&&&&\cr
&42&&30--35&&65--80&&&&&&&&&&&\cr
&43&&33--36&&&&&&&&&&&&&\cr
&44&&32--37&&68--83&&121--153&&&&&&&&&\cr
&45&&32--37&&80--84&&144--156&&242--268&&&&&&&\cr
&46&&34--38&&&&129--158&&&&&&&&&\cr
&47&&36--38&&&&116--161&&&&&&&&&\cr
&48&&34--39&&&&117--164&&&&&&&&&\cr
&49&&36--40&&66--90&&&&&&&&&&&\cr
&50&&40&&91--92&& && && && && &\cr
\vhop
}
\hrule
}
$$
\vfill\eject
%
%
%
%
%
%
%
%
\def\quad{\hskip 0.6em\relax}
\def\quod{\hskip 0.6em\relax}
\def\vhop{height2pt&\omit&&\omit&&\omit&&\omit&&\omit&\cr}
\noindent{\bf Table p=3.}
\smallskip
$$
\vcenter{
\tablefont
\lineskip=1pt
\baselineskip=10pt
\lineskiplimit=0pt
\setbox\strutbox=\hbox{\vrule height .7\baselineskip
                                depth .3\baselineskip width0pt}%
\offinterlineskip
\hrule
\halign{&\vrule#&\strut\quod\hfil#\quad\cr
height2pt&\omit&&\omit&&\omit&&\omit&&\omit&\cr
&$g\backslash q$&&3&&9&&27&&81&\cr
height2pt&\omit&&\omit&&\omit&&\omit&&\omit&\cr
\vhop
\noalign{\hrule}
\vhop
height2pt&\omit&&\omit&&\omit&&\omit&&\omit&\cr
&1&&7&&16&&38&&100&\cr
&2&&8&&20&&48&&118&\cr
&3&&10&&28&&55--58&&136&\cr
&4&&12&&28--30&&64--68&&154&\cr
&5&&12--14&&32--36&&55--78&&152--172&\cr
&6&&14--15&&29--40&&76--88&&190&\cr
&7&&16--17&&&&&&160--208&\cr
&8&&15--18&&38--47&&&&218--226&\cr
&9&&19&&40--51&&82--118&&244&\cr
&10&&19--21&&46--55&&&&226--262&\cr
\noalign{\hrule}
&11&&20--22&&55--59&&&&&\cr
&12&&22--24&&55--63&&109--148&&298&\cr
&13&&24--25&&60--66&&136--156&&&\cr
&14&&24--26&&56--70&&&&&\cr
&15&&28&&58--74&&&&&\cr
&16&&27--29&&74--78&&&&370&\cr
&17&&&&&&&&&\cr
&18&&24--31&&74-85&&&&&\cr
&19&&&&76--88&&&&&\cr
&20&&30--34&&74--91&&&&&\cr
\noalign{\hrule}
&21&&28--35&&82--95&&&&&\cr
&22&&28--36&&76--98&&&&&\cr
&23&&&&&&&&&\cr
&24&&28--38&&82--104&&&&&\cr
&25&&30--40&&&&&&&\cr
&26&&&&110--111&&&&&\cr
&27&&&&&&&&&\cr
&28&&36--43&&&&&&&\cr
&29&&&&&&&&&\cr
&30&&&&&&&&&\cr
\noalign{\hrule}
&31&&&&&&&&&\cr
&32&&36--48&&&&&&&\cr
&33&&&&&&&&&\cr
&34&&&&&&&&&\cr
&35&&&&&&&&&\cr
&36&&&&&&244--319&&730&\cr
&37&&&&109--145&&&&&\cr
&38&&&&&&&&&\cr
&39&&&&109--152&&244--339&&&\cr
&40&&&&&&244--346&&&\cr
\noalign{\hrule}
&41&&&&&&&&&\cr
&42&&&&&&&&&\cr
&43&&&&&&&&&\cr
&44&&&&&&&&&\cr
&45&&&&&&&&&\cr
&46&&55--63&&&&&&&\cr
&47&&&&&&&&&\cr
&48&&55--66&&&&&&&\cr
&49&&&&&&&&&\cr
&50&&&&182--186&&&&&\cr
\vhop
}
\hrule
}
$$
\vfill\eject
\vfill\eject
\bigskip
\noindent
\centerline{\bf References}
\bigskip
\noindent
\ni
[D-V] V.G.\ Drinfeld, S.G.\  Vladut : The number of points on an
algebraic curve. {\sl Funct. Anal. and Appl.\ \bf  17}
(1983),  p.\ 53-54.

\smallskip \ni
[G-S1]  A.\ Garcia , H. Stichtenoth: A tower of Artin-Schreier
extensions of function fields attaining the Drinfeld-Vladut bound.
{\sl Inv.\ Math.\ \bf 121} (1995), p.\ 211-222.

\smallskip \ni 
[G-S2]  A.\ Garcia , H. Stichtenoth: Algebraic function fields
over finite fields with many places. 
{\sl IEEE Trans.\ Inf.\ Th.\ \bf 41} (1995), p.\ 1547-1563.

\smallskip \ni
[G-V 1] G.\ van der Geer, M.\ van der Vlugt: Artin-Schreier curves
and codes. {\sl J.\ Algebra \bf 139} (1991),p.\ 256-272.

\smallskip \ni
[G-V 2] G.\ van der Geer, M.\ van der Vlugt : Reed-Muller
codes and supersingular curves. I. {\sl Compositio Math.\
\bf 84} (1992), p.\ 333-367.

\smallskip \ni
[G-V 3] G.\ van der Geer, M.\ van der Vlugt : Curves over finite
fields of characteristic $2$ with many rational points. {\sl Comptes
Rendus Acad.\ Sci.\ Paris \bf 317}, S\'erie I (1993), p.\ 593-597.

\smallskip \ni
[G-V 4] G.\ van der Geer, M.\ van der Vlugt : Generalized
Hamming weights of Melas codes and dual Melas codes.
{\sl SIAM Journal on
Discrete Math. \bf 7 } (1994), p.\ 554-559.

\smallskip \ni
[G-V 5] G.\ van der Geer, M.\ van der Vlugt : The second
generalized Hamming weight of $BCH(2)^{\bot}$.
{\sl Bull.\ London Math.\ Soc.\ \bf 27} (1995), p.\ 82-86.

\smallskip \ni
[G-V 6] G.\ van der Geer, M.\ van der Vlugt : Fibre
products of Artin-Schreier curves and generalized
Hamming weights of codes. {\sl J.
Comb. Theory  A \bf 70}(1995),p.\ 337-348.

\smallskip \ni
[G-V 7]  G.\ van der Geer, M.\ van der Vlugt : Generalized Hamming
weights of codes and curves over finite fields with many points.
Israel Math Conf. Proc. (to appear).

\smallskip \ni
[G-V 8]  G.\ van der Geer, M.\ van der Vlugt : Quadratic forms,
generalized Hamming weights of codes and curves with many points.
Preprint November 1994. To appear in J.\ of Number Theory.

\smallskip \ni
[G] V.D.\ Goppa : Codes on algebraic curves. {\sl Sov.\ Math.\  Dokl.\
\bf 24} (1981),p.\ 170-172.

\smallskip \ni
[I] Y.\ Ihara: Some remarks on the number of rational points of
algebraic curves over finite fields. {\sl J.\ Fac.\ Sci.\ Tokyo \bf
28} (1981), p.\ 721-724.

\smallskip \ni
[M] Y. Manin: What is the maximum number of points on a curve over
$\F_2$ ?  {\sl J.Fac.\ Sci.\ Tokyo \bf 28}, (1981) p.\ 715-720.

\smallskip \ni 
[N-X] H.\ Niederreiter, C.\ Xing: Cyclotomic function fields and explicit
global function fields with many rational places. Report Institute of
Information Processing, Austrian Acad. of Sciences, Vienna, 1995.

\smallskip \ni
[Sch] R.\ Schoof: Algebraic curves and coding theory. UTM 336, Univ.\
of Trento, 1990.

\smallskip \ni
[S 1] J-P.\ Serre : Sur le nombre de points rationnels d'une  courbe
alg\'ebrique sur un corps fini.  {\sl C.R.\ Acad.\ Sci.\ Paris \bf 296}
S\'erie I (1983), p. 397-402. (= Oeuvres III, p.\ 658-663).

\smallskip \ni
[S 2] J-P.\ Serre : Nombre de points des courbes
alg\'ebriques sur $\F_q$. S\'em.\ de Th\'eorie des
Nombres de Bordeaux, 1982/83, exp.\ no.\ 22. (= Oeuvres III, No.\
129, p.\ 664-668).

\smallskip \ni
[S 3] J-P.\ Serre : Quel est le nombre maximum de points
rationnels que peut avoir une courbe alg\'ebrique de genre
$g$ sur un corps fini $\F_q$ ? R\'esum\'e des Cours de 1983-1984.
(Oeuvres III, No.\ 132, p.\ 701-705).

\smallskip \noindent
[S-X] H.\ Stichtenoth, C. Xing: The genus of maximal function
fields. {\sl Manuscripta Math. \bf 86} (1995), p.\ 217-224.

\smallskip \ni
[T-V-Z]  M.A.\ Tsfasman, S.G.\ Vladut, T.\ Zink: Modular curves,
Shimura curves and Goppa codes better than the Varshamov-Gilbert
bound. {\sl Math.\ Nachr.\ \bf 109} (1982), p.\ 21-28.

\smallskip \ni
[W] V.K.\ Wei:
Generalized Hamming weights for linear codes. {\sl IEEE Trans.\
Inform.\ Th.\ \bf 37} (1991), p.\ 1412-1418.

\smallskip \ni
[Wi] M.\ Wirtz : Konstruktion und Tabellen linearer Codes.
Westf\"alische Wilhelms-Universit\"at M\"unster, 1991.

\smallskip \ni 
[X] C.\ Xing: Curves with many rational points over finite fields and Drinfeld
modules. Report Univ. of Science and Technology of China, Hefei, 1995.

\bigskip
\noindent
\bigskip
\settabs3 \columns
\+G. van der Geer  &&M. van der Vlugt\cr
\+Faculteit
Wiskunde en Informatica &&Mathematisch Instituut\cr
\+Universiteit van
Amsterdam &&Rijksuniversiteit te Leiden \cr
\+Plantage Muidergracht 24&&Niels Bohrweg 1 \cr
\+1018 TV Amsterdam
&&2300 RA Leiden \cr
\+The Netherlands &&The Netherlands \cr
\bigskip
\noindent
updated version:  May 14, 1996
\bye